\documentclass[review,times]{elsarticle}
\usepackage[a4paper,left=1in,right=1in,top=1in,bottom=1in,footskip=.25in]{geometry}

\usepackage{float} 
\usepackage{graphicx,multirow,soul}
\usepackage{amsmath}
\usepackage{epsfig}
\usepackage{xcolor}
\usepackage{amsfonts}
\usepackage{bm}
\usepackage[utf8]{inputenc}
\usepackage{color}
\usepackage{siunitx}
\usepackage{bm}
\biboptions{sort&compress}

\newcommand{\beq}{\begin{equation}}
\newcommand{\eeq}{\end{equation}}
 \newcommand{\eps}{\varepsilon}

\newcommand{\dis}{\displaystyle}
\newcommand{\RR}{{\rm I\kern - 1.5pt R}}

\def\tb{\textbf}

\def\tb{\textbf}
\newcommand{\bs}[1]{\boldsymbol{#1}}

\journal{Journal of Computational and Applied Mathematics}

\begin{document}
\begin{frontmatter}
\title{Sparse Polynomial Chaos Expansion for Universal Stochastic Kriging}
\author[1] {J. C. Garc\'{\i}a-Merino}
\author[1] {C. Calvo-Jurado}
\author[2] {E. Garc\'{\i}a-Mac\'{\i}as}
\address[1]{Department of Mathematics, School of Technology, 10003, C\'aceres, Spain}
\address[2]{Department of Structural Mechanics and Hydraulic Engineering, University of Granada, Campus de Fuentenueva s/n, 18071 Granada, Spain}

\begin{abstract}
Surrogate modelling techniques have opened up new possibilities to overcome the limitations of computationally intensive numerical models in various areas of engineering and science. However, while fundamental in many engineering applications and decision-making, the incorporation of uncertainty quantification into meta-models remains a challenging open area of research. To address this issue, this paper presents a novel stochastic simulation approach combining sparse polynomial chaos expansion (PCE) and Stochastic Kriging (SK). Specifically, the proposed approach adopts adaptive sparse PCE as the trend model in SK, achieving both global and local prediction capabilities and maximizing the role of the stochastic term to conduct uncertainty quantification. To maximize the generalization and computational efficiency of the meta-model, the Least Angle Regression (LAR) algorithm is adopted to automatically select the optimal polynomial basis in the PCE. The computational effectiveness and accuracy of the proposed approach are appraised through a comprehensive set of case studies and different quality metrics. The presented numerical results and discussion demonstrate the superior performance of the proposed approach compared to the classical ordinary SK model, offering high flexibility for the characterization of both extrinsic and intrinsic uncertainty for a wide variety of problems. 
\end{abstract}


\end{frontmatter}


\section{Introduction}


Computer models, commonly referred to as simulators, have become ubiquitous assets in most areas of engineering and science. In cases where systems or processes are assumed to be predictable and devoid of inherent randomness, a deterministic simulator can be adopted. Nonetheless, deterministic approaches appear inadequate for facing real-world problems characterized by intrinsic or non-parameterisable randomness, for which the adoption of stochastic simulators becomes necessary for proper modelling. A stochastic simulator incorporates randomness or uncertainty in its predictions. Thus, while a deterministic simulator always produces the same prediction for a fixed input, in stochastic simulators, each simulation run may produce different results as a result of the inherent randomness in the model. Such models have gained significant popularity with manifold applications in diverse fields such as epidemiology~\cite{Allen-2017, Mroue-2021}, materials science~\cite{Thiedmann-2011,Meyer-2003}, or reliability analysis~\cite{Sakki-2022}, enabling researchers to make robust predictions to assist decision-making for complex systems. However, for both deterministic and stochastic computationally demanding simulators, formidable challenges arise when implemented in iterative procedures such as optimization~\cite{Luo-2011}, sensitivity analysis~\cite{Khaledi-2016}, parameter inference~\cite{Garcia-2020}, or uncertainty propagation~\cite{Gar-2022}. Despite the introduction of numerous surrogate models or meta-models in the literature over the past few years to replace intensive numerical models and alleviate the computational burden, their success mostly limits to deterministic simulators, while the consideration of stochasticity remains an open research topic.

There is a vast literature on the development of surrogate models for deterministic simulators~\cite{Alizadeh-2020, Kudela-2022}. In particular, non-intrusive models such as polynomial chaos expansion (PCE)~\cite{Bla-2011}, radial basis functions~\cite{Nunez-2018} or Kriging~\cite{Kleijnen-2009} have received significant attention from the scientific community due to their effectiveness in addressing resource-intensive simulations. The construction of these meta-models is usually conducted on the basis of a training dataset, also referred to as the experimental design (ED)~\cite{Kleijnen-2009}. Such a training dataset contains the evaluations of the forward model, representing the quantity of interest (QoI), for a set of combinations of the input parameters covering the parameter space. In the context of stochastic simulators, the QoI usually presents a certain level of uncertainty, stemming from the stochastic nature of the model parameters and/or the presence of epistemic uncertainties in the model predictions~\cite{Baker-2022}. A first approach to meta-model such processes involves estimating the probability density function (PDF) of the simulator's predictions through functional decomposition~\cite{Moutoussamy-2015}. Alternatively, Stochastic Kriging (SK)~\cite{Ank-2010} has gained significant recognition as a general meta-modelling tool for representing stochastic simulation response surfaces. Defined as an extension of Kriging, SK accounts for heteroscedastic noise in the response by estimating the sample variance at each point of the ED using multiple runs of the simulator for each specific input. The effectiveness of SK has been demonstrated on a broad range of applications, including risk and reliability analysis~\cite{Chen-2016, Hao-2021}, industrial welding processes~\cite{Ruan-2018}, and game theoretic simulations~\cite{Pousi-2010}, to mention a few.

One of the primary challenges encountered in SK is the requirement of a significant number of replications at each ED point to accurately estimate the sample variance~\cite{Baker-2022}, which often leads to serious computational demands. One approach to alleviate such a limitation regards the optimal definition of the ED. In the case of stochastic simulators, sampling efforts are divided into two parts:     (i) sampling new ED locations, and  (ii) acquiring additional replications at existing points. In this light, numerous sequential sampling techniques have been presented   (refer to~\cite{Chen-2017, Wang-2018} for an extensive state-of-the-art review). On the other hand, various extensions of the classical SK have been proposed in the literature. It is worth noting the work by Binois \textit{et al.}~\cite{Binois-2018}, who proposed modelling the input-dependent noise as a separate Gaussian process. Additionally, they employed computationally efficient techniques such as the Woodbury identity to address the challenges encountered in SK for large EDs. In a similar vein, Hao \textit{et al.}~\cite{Hao-2021} introduced a Nested SK model, which decouples in an iterative way the response noise parameters from the deterministic Kriging parameters used for calibrating the meta-model. In a similar context, Cheng \textit{et al.}~\cite{Che-2021} proposed a generalized polynomial chaos-stochastic Kriging (gPCE SK) approach, which combines the benefits of PCE and SK. The gPCE SK method proposed by those authors utilizes a local PCE surrogate model at each design point, which is calibrated using a small number of simulation runs. By relying on these local surrogate models, numerous replications can be performed at each training point, significantly reducing the computational burden involved in the calibration of the SK model. However, it is important to note that this model assumes a distinction between controllable design variables and non-controllable noise variables, the latter causing the variations between replicates at each fixed design point. Furthermore, it relies on knowledge of the PDFs of the uncontrollable variables, which may be a strong assumption in many real-life applications.

In light of the discussion above, it is apparent that the development and application of surrogate modelling approaches for stochastic simulators have experienced considerable advances in recent years. Nevertheless, there are some aspects of SK that remain not fully addressed in the literature. Notably, while most applications of SK in the literature are limited to the use of constant trend models (Ordinary Kriging), the use of more advanced models for Universal SK remains largely unexplored in the literature, to the best of the author's knowledge. To address this gap, this paper presents an extension of the deterministic Polynomial-Chaos Kriging (PCK) method originally presented in~\cite{Schobi-2015} for the meta-modelling of  stochastic simulators with predictions contaminated by noise. The PCK meta-model is a powerful technique that combines PCE with Kriging to develop highly efficient surrogate models for complex systems. Notably, a prior work by the authors~\cite{Gar-2023} demonstrated the ability of PCK meta-modelling to successfully capture the behaviour of highly non-linear functions. Drawing from this positive outcome, the current study aims to apply the principles of PCK to the realm of stochastic simulators. Specifically, the proposed meta-model involves a Universal SK model with an adaptive sparse PCE model as the trend term. The orthonormal polynomial basis of the PCE is optimally selected by a model selection technique for sparse linear models, the least-angle regression (LAR) algorithm set out by Efron~\cite{Efron-2004}. Afterwards, the adjusted LAR-PC model is inserted as the trend term in a Universal SK model, and the hyper-parameters of the resulting meta-model are calibrated through a genetic algorithm (GA) global optimization approach. The effectiveness of the proposed approach is validated  through a set of well-known benchmark case studies in the surrogate modelling literature. The presented results and discussion demonstrate the potential of the proposed approach and the importance of   selecting a suitable trend model in the context of SK. 

The remainder of this paper is organized as follows. Section~\ref{K_and_SK_sec} outlines the theoretical formulation of Kriging and SK respectively. Section~\ref{Sect_3} provides a concise review of PCE and LAR and illustrates the proposed LAR-PCE SK approach. Section~\ref{Sect_4} presents the numerical results and discussion and, finally, Section~\ref{Sect_5} closes the paper with concluding remarks and potential future developments.

\section{Surrogate modelling: Deterministic and Stochastic Kriging} \label{K_and_SK_sec}

Let $(\Omega,{\cal F},\mu)$ be a probability space, where $\Omega$ denotes the event space equipped with a $\sigma$-algebra $\mathcal{F}$ of subsets of $\Omega$ and the probability measure $\mu$. Given $\Omega\subset \mathbb{R}^M$ as a connected set, let ${\mathcal M}:\Omega \subset \mathbb{R}^M\to\mathbb{R}$ be a computational model mapping an $M$-dimensional vector of input (design) variables $\bs{x} = \left[x_1,\ldots,x_M\right]^\textrm{T}$ into the output variable $ y\in\mathbb{R}$ (response surface or  QoI). The aim of a surrogate model is to replace the original model $\mathcal{M}$ by an approximated one $\mathcal{M}^*$, whose evaluation needs lower computational cost. When defined in a non-intrusive framework, the hyper-parameters of ${\cal M}^*$ must be calibrated on a basis of a certain training dataset. Such a training dataset consists in the realization of the forward model $\cal M$ over a set ED=$\left\{ \bs{x}^{(1)}, \ldots , \bs{x}^{(k)} \right\}$ of $k$ samples of the input variable $\bs{x}$ covering the design space, often referred to as the {\em experimental design} (ED). Note that, when the model is deterministic, the purpose of the surrogate model is simply to approximate the response surface $y = \mathcal{M}({\bs{x}})$ $\forall \bs{x} \in \Omega$. Conversely, when the model $\mathcal{M}$ exhibits a stochastic behaviour and its response is affected by a zero-mean random noise, the goal is to approximate the mean surface $\mathbb{E}[\mathcal{M}(\Omega)]$  and, possibly, provide a certain characterization of the uncertainty ~\cite{Ank-2010}. In this context, Kriging or \emph{Gaussian process} meta-modelling has proved to represent a powerful meta-modelling approach commonly used in engineering or geostatistics. In a deterministic computer experiment, the Kriging model assumes that the model response ${\mathcal M}(\bs{x})$ can be conceived as the sum of a deterministic and a stochastic model as:  

\begin{equation}\label{Krig1}
\mathcal{M}\left( \bs{x} \right) \approx \mathcal{M}^K \left( \bs{x} \right)=\mathcal{T}\left( \bs{x} \right)+\mathcal{Z}\left( \bs{x} \right),
\end{equation}

\noindent where $\mathcal{Z} \left( \bs{x} \right)$ is a zero-mean stochastic process and the trend $\mathcal{T}\left( \bs{x} \right) = \bs{f} \left( \bs{x} \right)^\textrm{T} \bs{\beta} = \dis\sum_{r=1}^{s} \beta_r f_r \left( \bs{x} \right)$ is a linear regression model with $s$ regression coefficients $\beta_r$ and user-selected regression functions $f_r \left( \bs{x} \right)$ (e.g $f_r \left( \bs{x} \right) = \bs{x}^r$). The formulation in Eq.~(\ref{Krig1}) corresponds to the most general version of Kriging, often referred to as {\em Universal Kriging} (UK). Simpler expressions known as {\em Simple Kriging} (SK) and {\em Ordinary Kriging} (OK) are deduced when the trend term is assumed to be constant throughout the design space. In SK, it is assumed that the mean of the variable being predicted is known and constant across the entire  area of study, whereas (OK) estimates  the mean based on the available sample data. An interpretation of Eq.~(\ref{Krig1}) is that deviations from the regression model may resemble a sample path of a properly chosen stochastic process~\cite{DACEbasicos}. The stochastic process $\mathcal{Z} \left( \bs{x} \right)$ is determined by the process variance $\sigma$ and an auto-correlation function $\tb{R} \left( \bs{x},\,\bs{x}' \right) = \mathbf{R}\left( \left|{\bs{x}-\bs{x}'} \right| ;\bs{\theta} \right)$, which depends on the spatial distance $\left| \bs{x}^{(i)}-\bs{x}^{(j)} \right|$ and some hyper-parameters $\bs{\theta}$ to be determined. The auto-correlation function $\textbf{R}$ must satisfy ${\bf{R}} \left(\bs{0};\bs{\theta} \right) = 1$ and ${\bf{R}} \left(\bs{x}^{(i)}-\bs{x}^{(j)};\bs{\theta} \right)\to 0$ when $\left|\bs{x}^{(i)}-\bs{x}^{(j)} \right| \to +\infty$. Thus, ${\cal Z}$ is characterized by the covariance matrix~\cite{San-2003, Ras-2006}:

\begin{equation}\label{esti_Z}
\bs{\Sigma}_\mathcal{Z} :=  \left\{{\rm Cov} \left({\mathcal{Z}} \left(\bs{x}^{(i)} \right),{\mathcal{Z}} \left(\bs{x}^{(j)}\right)\right)\right\}_{ij}={\sigma}^2 \, \bf{R} \left( \left|\bs{x}^{(i)}-\bs{x}^{(j)}\right|; \bs{\theta} \right).
\end{equation}

Among the possible choices of the function $\textbf{R}$, a common alternative and the one used in the present work is the Gaussian correlation function $\mathbf{R}\left( \left|{\bs{x}-\bs{x}'} \right|; \bs{\theta} \right) = \rm{exp} \left(-\bs{\theta} \left|\bs{x}-\bs{x}'\right|^2 \right)$~\cite{DACEbasicos}.

The Kriging model is thus fully determined by the hyper-parameters $\left\{\bs{\beta}, \bs{\theta},\sigma \right\}$, which are usually estimated by maximizing the likelihood function over the $k$ samples of the  ED~\cite{marrel}:

\begin{equation}\label{lk}
\mathcal{L} \left( \bs{\theta},\sigma,{\bs{\beta}}; \bf{y} \right) = \frac{1}{\sqrt{ \left( 2\pi\sigma^2 \right)^k \left| \textbf{R} \right| }} \, \rm{exp} \left[-\frac{1}{2\sigma^2} \left( \textbf{y}-\textbf{F}\bs{\beta}\right)^\textrm{T} \textbf{R}^{-1} \left(\textbf{y}-\textbf{F}\bs{\beta} \right)
\right],
\end{equation} 

\noindent where ${\tb y} = \left[ {y}^{(1)}={\mathcal M}\left(\bs{x}^{(1)}\right),\ldots, {y}^{(k)}={\mathcal  M}\left(\bs{x}^{(k)}\right) \right]^{\rm T}$ and $\textbf{F}_{ij} = f_j(\bs{x}^{(i)})$ is referred to as the {\em information} matrix. Nevertheless, it is possible to extract estimates in  closed form of the hyper-parameters $\bs{\hat{\beta}}$ and $\bs{\hat{\sigma}}$ as functions of the estimates of $\hat{\bs{\theta}}$ by the best linear unbiased estimator (BLUE) as:

\begin{equation}\label{calibrar}
\bs{\hat{\beta}}(\hat{\bs{\theta}}) = \left( \textbf{F}^\textrm{T}\textbf{R}^{-1}\textbf{F}\right)^{-1}\textbf{F}^\textrm{T}\textbf{R}^{-1}\textbf{y} ,\quad \hat{\sigma}^2 \left(\hat{\bs{\theta}}\right) = \frac{1}{k} \left( \textbf{y}-\textbf{F}{\bs{\beta}} \right)^{\textrm{T}} \textbf{R}^{-1} \left( \textbf{y} -\textbf{F}{\bs{\beta}} \right). 
\end{equation} 

Substituting Eq.~(\ref{calibrar}) into (\ref{lk}) and taking logarithm, a simplified optimisation problem that depends only on ${\bs \theta}$ is derived as:

\begin{equation}\label{MLCV}
\hat{\bs{\theta}} = \underset{\bs{\theta}}{\textrm{arg\,min}} \left[\frac{1}{k} \left(\textbf{y} -\textbf{F}{\bs{\beta}} \right)^{\textrm{T}} \textbf{R}^{-1} \left( \textbf{y} -\textbf{F}{\bs{\beta}} \right) \left| \textbf{R} \right|^{1/k} \right].
\end{equation} 

Therefore, the Kriging mean and mean-squared error (MSE) at any point $\bs{x}$ in the parameter space can be computed as~\cite{Ras-2006}:

\begin{equation}
\label{kriging_pred}
\begin{array}{c}
\hat{\cal M}^K(\bs{x}) = \bs{f}(\bs{x})^\textrm{T}\hat{\bs{\beta}}+\bs{r}(\bs{x})^\textrm{T}\textbf{R}^{-1} \left({\bs y}-\textbf{F}\hat{\bs{\beta}}\right),\\
\widehat{\rm{MSE}}(\bs{x})=\sigma^2 \left[ 1-\bs{r}(\bs{x})^\textrm{T}\textbf{R}^{-1}\bs{r}(\bs{x})+\bs{u}(\bs{x})^\textrm{T}\left(\textbf{F}^\textrm{T}\textbf{R}^{-1}\textbf{F}\right)^{-1}\bs{u}(\bs{x})
\right],
\end{array}
\end{equation} 

\noindent where $\bs{u} \left( \bs{x} \right) = \textbf{F}^\textrm{T}\textbf{R}^{-1}\bs{r}(\bs{x})-\bs{f}(\bs{x})$ and $r_i(\bs{x}) = \textbf{R} \left( \left| \bs{x}-\bs{x}^{(i)} \right|; \hat{\bs{\theta}} \right)$, $i = 1, \ldots, k$.

 To accommodate the existence of stochasticity in the response surface model ${\cal M}(x)$, Ankenman \textit{et al.}~\cite{Ank-2010} proposed Stochastic Kriging as   an extension of the classical Kriging scheme specified in Eq.~(\ref{Krig1})  to define the response surface $\mathcal{M}(\bs{x})$ in stochastic terms. In this context, several simulation runs in the same design point are performed, so the experimental design must be redefined as $\textrm{ED} = \left\{ \left( \bs{x}^{(i)}, n_i \right) \right\}_{i=1}^{k}$, where $n_i \in \mathbb{N}$ is the number of replications at the design point $\bs{x}^{(i)}$. In particular, the output of the $l$-th simulation replication at $\bs{x}^{(i)}$, $i=1,\ldots,k$ is given by:
\begin{equation}\label{sk}
\mathcal{M}_i \left( \bs{x}^{(i)} \right) ={\bs f}(\bs{x})^\textrm{T}{\bs{\beta}}+\mathcal{Z}\left(\bs{x}^{(i)}\right)+\eps_l \left( \bs{x}^{(i)} \right),\quad l = 1, \ldots, n_i, 
\end{equation}
\noindent where the terms $\mathcal{Z}\left(\bs{x}^{(i)}\right)$ and $\eps_l \left( \bs{x}^{(i)} \right)$ represent the model {\em extrinsic uncertainty} and the {\em intrinsic uncertainty} for the $l$-th replication, respectively. In other words, the changes in the model response are due to the simulator evaluation at different points in the parameter space, or due to the stochastic nature of the model itself, respectively. Analogously to Eq.~(\ref{esti_Z}), the intrinsic uncertainty can be characterized by the $k\times k$ covariance matrix:
\begin{equation}\label{Sigma_eps_defi}
\bs{\Sigma}_{\eps} = \left\{ {\rm{Cov}} \left(\bar{\eps}\left({\bs{x}}^{(i)}\right),\bar{\eps}\left({\bs{x}}^{(j)}\right) \right) \right\}_{ij}, \quad \overline{\eps}\left(\bs{x}^{(i)}\right)=\dfrac{1}{n_i} \dis\sum_{l=1}^{n_i}{\eps_l}\left(\bs{x}^{(i)}\right), \, i=1,\ldots,k.
\end{equation}
In this light, similarly to Kriging, the calibration of the surrogate model involves maximizing a likelihood function, whose expression is the natural consequence of incorporating the intrinsic uncertainty in the Kriging methodology: 
\begin{equation}\label{l-SK}
\mathcal{L} \left( \bs{\beta},\sigma^2,\bs{\theta};\overline{\cal M} \right) = \frac{1}{\sqrt{\left|\bs{\Sigma}_{\cal Z}+\bs{\Sigma}_{\eps}\right| \, (2\pi)^{k}}} \, {\rm exp} \left(-\frac{1}{2} \left(\overline{\cal M}-\textbf{F}\bs{\beta} \right)^\textrm{T} \left(\bs{\Sigma}_{\cal Z}+\bs{\Sigma}_\eps\right)^{-1} \left(\overline{\cal M}-\textbf{F}\bs{\beta} \right) \right),
\end{equation}
\noindent where 
\begin{equation}
\overline{\mathcal{M}} = \left[\overline{\mathcal{M}}\left(\bs{x}^{(1)}\right) ,\ldots, \overline{\mathcal{M}}\left(\bs{x}^{(k)}\right) \right]^\textrm{T}, \quad \overline{\mathcal{M}}\left(\bs{x}^{(i)}\right) = \frac{1}{n_i} \dis\sum_{l=1}^{n_i}\mathcal{M}_l(\bs{x}^{(i)}),
\end{equation}

\noindent and $\bs{\Sigma}_{\mathcal Z}(\bs{x},\cdot)$ represents the $k\times 1$ {\em cross-correlation vector} $\left\{ {\rm{Cov}} \left(\mathcal{Z} \left(\bs{x} \right),\mathcal{Z} \left(\bs{x}^{(j)} \right)\right) \right\}_{j=1}^k$ between the prediction point $\bs{x}$ and each one of the ED points. Upon the computation of the SK parameters, the SK predictor along with its associated mean squared error (MSE) can be derived as follows~\cite{Chen-2014}: 
\begin{equation}\label{sk-pred}
\begin{array}{c}
\hat{\mathcal{M}}^{SK}\left( \bs{x} \right)={\bs f}\left( \bs{x} \right)^\textrm{T}\hat{\bs{\beta}} +\bs{\Sigma}_{\mathcal Z}\left( \bs{x},\cdot \right) \left[ \bs{\Sigma}_{\mathcal Z}+\bs{\Sigma}_{\eps} \right]^{-1} \left( \overline{\mathcal{M}}-\textbf{F}\hat{\bs{\beta}} \right),
\\
\widehat{\rm MSE}\left(\bs{x}\right) = \bs{\Sigma}_{\mathcal
Z}\left(\bs{x},\bs{x}\right)- \bs{\Sigma}_{\mathcal
Z}\left(\bs{x},\cdot\right)^\textrm{T}\left[\bs{\Sigma}_{\mathcal
Z}+\bs{\Sigma}_{\eps}\right]^{-1} \bs{\Sigma}_{\mathcal
Z}\left(\bs{x},\cdot\right)+{\bs\zeta}^\textrm{T} \left( \textbf{F}^\textrm{T}  \left[\bs{\Sigma}_{\mathcal
Z}+\bs{\Sigma}_{\eps}\right]^{-1}\textbf{F} \right)^{-1}{\bs\zeta}, 
\end{array}
\end{equation}

\noindent with $\bs{\zeta}= {\bs f} \left( \bs{x} \right) - \textbf{F}^\textrm{T} \left[ \bs{\Sigma}_{\mathcal Z}+\bs{\Sigma}_{\eps} \right]^{-1} {\bs{\Sigma}_{\cal Z} \left( \bs{x},\cdot \right)}$. It is worth mentioning that the  last term in Eq.~(\ref{sk-pred}) represents the variability inflation of the SK predictor resulting from the estimation of $\hat{\bs{\beta}}$. Additionally, it is important to highlight that, in deterministic scenarios where the intrinsic uncertainty is negligible, $\bs{\Sigma}_{\eps}$ vanishes and Eq. (\ref{sk-pred}) naturally converges to Eq. (\ref{kriging_pred}). In other words, SK reduces to the standard Kriging metamodel when $\bs{\Sigma}_{\eps}$ tends to zero.

\section{LAR-PCE based Stochastic Kriging} \label{Sect_3}

The higher complexity of UK compared to OK might suggest that it leads to improved estimates. However, as reported in the literature, a trend misspecification can lead  UK not to  improve the accuracy in its estimates~\cite{Taj-2015} or even  worsen its performance~\cite{Staum-2009}. This limitation can be attributed to the fact that UK involves the estimation of additional parameters ($\bm{\beta}$), which can potentially induce over-fitting~\cite{Kle-2016}. Therefore, to harness the advantages of UK while mitigating its potential drawbacks, it is crucial to choose a minimal trend basis that effectively captures the complexity of the response surface. Motivated by these considerations, the objective of this section is to identify an optimal choice of the trend model $\bs{f}$ in Eq.~(\ref{sk}) by combining PCE and LAR. Hereafter, the fundamentals of PCE and the LAR algorithm are briefly presented, followed by their combination into Universal SK.

\subsection{Polynomial-chaos expansion}\label{PCE_sec}

The PCE model represents the model output response $y\in \mathbb{R}$ as its expansion through a basis of orthonormal multivariate polynomials $\Psi_{\bm\alpha}$ as:
\begin{equation}\label{expansion}
y = \mathcal{M} \left( \bs{x} \right) = \dis\sum_{\bm{\alpha} \, \in \, \mathbb{N}^M} a_{\bm{\alpha}} \Psi_{\bm{\alpha}}(\bs{x}),
\end{equation} 
\noindent where $a_{\bm{\alpha}}$ are the coefficients of the expansion, and the multi-dimensional index notation $\bm{\alpha} = \left[\alpha_1,\ldots,\alpha_M \right]^\textrm{T}$, $\alpha_i \in \mathbb{N}$, has been adopted. Depending on the nature of the input variables, classical families of orthonormal polynomials $\Psi_{\bm\alpha}$ include the Hermite or Legendre polynomials~\cite{Xiu-2002}. Although the expression~(\ref{expansion}) is exact on deterministic simulators for an infinite number of terms, only a finite number of terms can be computed in practice. Consequently, different truncation schemes have been proposed in the literature. One of the simplest approaches consists in selecting the set of all polynomials whose total degree  $\left| \bs{\alpha} \right| = \dis\sum_{i=1}^M \alpha_i$ is up to $p$, that is $\mathcal{A}^{M,p}=\{\bs{\alpha} \in \mathbb{N}^M:\,0 \leq \left| \bs{\alpha} \right| \leq p\}$ ~\cite{Faj-2017}. This leads to a total number of ${\cal P}$ terms in the expansion:
\begin{equation}\label{stability}
{\cal P}: = \rm{card}\left(\mathcal{A}^{M,p}\right) = \left( \begin{matrix}
M+p\\
p\\
\end{matrix}
 \right).
\end{equation}
Note however that, for high-dimensional or non-linear problems, this truncation scheme leads to computing a large number of coefficients. As a solution, a more restrictive {\em hyperbolic truncation scheme} was proposed by Blatman and Sudret~\cite{Bla-2011}. Based on the common observation in practice that coefficients associated with high-order interaction terms are usually close to zero (often referred to as the \emph{sparsity-of-effect principle}), this truncation scheme selects all the multi-indices belonging to the set $\mathcal{A}^{M,p,q} = \{\bs{\alpha} \in \mathbb{N}:\,0\leq \left\| \bs{\alpha} \right\|_q\leq p\}$ with $q$-norm $\left\|\bs{\alpha}\right\|_q = \left(\dis\sum_{i=1}^M \left| \alpha_i \right|^q\right)^{\frac{1}{q}}$. Note that, as the $q$-norm decreases, the number of mixed-order terms in the expansion also decreases until,   ultimately leaving only the univariate polynomials remain. 

Once the set of candidate polynomials and an ED of $k$ samples are defined, the expansion coefficients ${\bm a} = \left\{{\bm a_\alpha},\; {\bm \alpha} \in \mathcal{A}^{M,p} \subset \mathbb{N}^M \right\}$ are obtained by minimizing the expectation of the least squares  residuals:

\begin{equation}\label{leastsqdis} 
\bm{a}= {\textrm{arg} \min_{\bm{a} \, \in \, \mathbb{R}^{\mathcal P}}} \, \frac{1}{k} \dis\sum_{i=1}^k\left[{\mathcal M}\left(\bs{x}^{(i)}\right)- \dis\sum_{{\bm\alpha} \, \in \, {\mathcal A}^{M,p}}{\bm a}_{\bm \alpha}{\bs\Psi}_{\bm\alpha} \left(\bs{x}^{(i)} \right)\right]^2.
\end{equation}
Collecting the realizations of the QoI in the ED in vector form as ${\tb y} = \left[{y}^{(1)}={\mathcal M}(\bs{x}^{(1)}),\ldots, {y}^{(k)}={\mathcal  M}(\bs{x}^{(k)}) \right]^{\rm T}$, the solution of the optimization problem~(\ref{leastsqdis}) is given by:
\begin{equation}\label{sol}
\begin{array}{c}
\hat{\bm{a}}=\left({\bs{\Psi}}^{\textrm{T}} {\bs{\Psi}} \right)^{-1}{\bs{\Psi}}^{\textrm{T}}{\bm y}, \quad {\bs{\Psi}} = \left( {\Psi}_{ij} \right) = \left[\Psi_{\bm{\alpha}_j}\left(\bs{x}^{(i)}\right) \right]_{\, i=1,\ldots,\,k}^{\, j=1,\ldots,\,{\mathcal  P}},
\end{array}
\end{equation}
\noindent where the {\em information matrix} $\bm{\Psi}$ contains the evaluation of the basis polynomials on the ED. Once Eq.~(\ref{leastsqdis}) has been solved, the estimates by the resulting PCE surrogate model read:
\begin{equation}\label{PCE} 
{\mathcal{M}}^{PC}(\bs{x}) = \dis\sum_{\bm{\alpha} \, \in \, \mathcal{A}^{M,p}} \hat{\bm{a}}_{\bm{\alpha}} {\bs\Psi}_{\bm{\alpha}}(\bs{x}).
\end{equation}
The quality of the resulting expansion is typically assessed through cross-validation. Among the different metrics available in the literature, the leave-one-out (LOO) error $\epsilon_{LOO}$ is particularly interesting since it can be estimated in closed form for linear models. In the case of the PCE model as given in Eq.~(\ref{PCE}), $\epsilon_{LOO}$ can be expressed as follows:

\begin{equation}\label{LOO-eq}    
\epsilon_{LOO} = \frac{\dis\sum_{i=1}^k\left(\frac{\mathcal{M} \left( \bs{x}^{(i)} \right)- \mathcal{M}^{PC}(\bs{x}^{(i)}) }{1-h_i}\right)^2}{\dis\sum_{i=1}^k \left( \mathcal{M} \left( \bs{x}^{(i)} \right) - \hat{\mu}_{\bs{y}} \right)^2},
\end{equation}

\noindent where $h_i$ is the $i$-th component of the vector $\mathbf{h}={\rm{diag}}\left\{\bs{\Psi}\left({\bs{\Psi}}^{\textrm{T}} {\bs{\Psi}} \right)^{-1}{\bs{\Psi}}^{\textrm{T}} \right\}$ and $\hat{\mu}_{\bs{y}} = \frac{1}{k} \dis\sum_{i=1}^k \mathcal{M} \left( \bs{x}^{(i)} \right)$.

\subsection{Least angle regression (LAR)}\label{sec_LAR}

The LAR algorithm is a model selection approach~\cite{Efron-2004} that provides a very efficient procedure to implement the {\em Least Absolute Shrinkage and Selection Operator} (LASSO) method introduced by Tibshirani~\cite{Tib-1996}. LASSO is a procedure for shrinking some regression coefficients to zero in a regression model by adding a restriction or penalty term to the least squares problem as:

\begin{equation}\label{lasso}
\min_{{\bs{\beta}} \, \in \, \RR^{\cal P}} \left(   \left\| {\bs{\Psi}} \left(\bs{x}\right) {\bs{\beta}}-{\cal M} \left( \bs{x} \right) \right\|_2+\eta \|{\bs{\beta}}\|_1 \right), \quad \eta \in \RR, \quad \eta \geq 0.
\end{equation}

The norm $\left\| \cdot \right\|_1$ adds a penalty equal to the absolute value of the magnitude of the coefficients. In other words, it limits the size of the coefficients leading to sparse models, i.e.~models with few coefficients. The computation of the solution of the optimisation problem in Eq.~(\ref{lasso}) is a convex quadratic programming problem with linear inequality constraints, which can be solved very efficiently by means of the LAR algorithm~\cite{Efron-2004}. It is based on a variable selection procedure called {\em Incremental Forward Stepwise Regression}~\cite{Has-2009} that builds the model sequentially, adding one basis polynomial function $\Psi_{\bs{\alpha}_i}$ at a  time to the set ${\cal A}_i$ of active basis functions, i.e. the set of predictors most correlated with the current residual $r_i=y-{\bs{\Psi}}_{{\cal A}_i} \, {\bs{\beta}}_{{\cal A}_i}$  at the $i$-th iteration. The LAR algorithm can be summarily described by the following steps:

\begin{enumerate}
\item Set the vector residuals ${r}_0={y}$ and initialize the coefficients vector ${\bs{\beta}}_{\bs{\alpha}_i}=0$, $i=1,\ldots,{\cal P}-1$. Set all the predictors $\Psi_{\bs{\alpha}_i}$ standardized to have mean zero and unit norm. 
\item Find the predictor $\Psi_{\bs{\alpha}_i}$ most correlated with the current residual $r_i=y-{\bs{\Psi}}_{{\cal A}_i} \, {\bs{\beta}}_{{\cal A}_i}$ (with  ${\bs{\beta}}_{{\cal A}_i}$ the vector of coefficients for the set of the active basis functions ${{\mathcal A}_i}$ at the $i$-th iteration).
\item Move ${\bs{\beta}}_{\bs{\alpha}_i}$  from $\textbf{0}$ towards its least-squares coefficient on $\Psi_{\bs{\alpha}_i}$, until another predictor $\Psi_{\bs{\alpha}_j}$ has as much correlation with the current residual $r_i$ as $\Psi_{\bs{\alpha}_i}$.
\item Move $\{{\bs{\beta}}_{\bs{\alpha}_i},\,{\bs{\beta}}_{\bs{\alpha}_j}\}$ in the direction defined by their joint least squares coefficient of the current residual on $\{\Psi_{\bs{\alpha}_i},\,\Psi_{\bs{\alpha}_j}\}$, until another predictor $\Psi_{\bs{\alpha}_l}$ has the same correlation with the current residual. 
\item Continue the process until  $m$= min$\, ({\cal P},k-1)$ predictors have been added.
\end{enumerate}

Steps 3 and 4  update the set of {\em active} coefficients toward their least-square value in the form ${\bs{\beta}}_{{\cal A}_{i+1}}={\bs{\beta}}_{{\cal A}_{i}}+\zeta_{i}\bm{u}_{{\cal A}_{i}}$. The factor $\zeta_{i}$ is the smallest positive scalar such that some new index $j$ joins the active set ${\cal A}_{i+1}={\cal A}_{i}\cup \{j\}$, and the vector $\bm{u}_{{\cal A}_{i}}$ is known as the {\em descent direction} (refer to~\cite{Efron-2004} for further details). The name {\em least angle} arises from a geometric interpretation of this process, whereby $\bm{u}_{{\cal A}_{i}}$ provides the smallest (and equal) angle with each of the predictors in ${\cal A}_{i}$. Using the LAR scheme, a collection of $m$ PCE meta-models of decreasing sparsity  is extracted. Then, the optimal meta-model can be selected by a cross-validation scheme, to which the LOO error from Eq.~(\ref{LOO-eq}) has been adopted in this work. Typically, the optimal model maximizing prediction accuracy and minimizing overfitting presents a cardinality, i.e. the number of active basis funcions, lower that the cardinality $\cal P$ of the whole candidate basis. It is known that LAR only requires ${\mathcal O}({\cal P}^3 + k{\cal P}^2)$ computations~\cite{Efron-2004}. Therefore, when ${\cal P}\ll k$, it follows that ${\cal P}^3<k {\cal P}^2$ and ${\mathcal O}(k{\cal P}^2)\sim {\mathcal O}(k)$. Thus, the computational complexity added by the PCE model when inserted as the trend term in the SK model is marginal with respect to  the optimization problem related to the determination of the hyper-parameters in Eq.~(\ref{MLCV}),  which is ${\mathcal O}(k^3)$.

\subsection{LAR-PCE Stochastic Kriging}\label{sec_LAR-PCE SK}

\begin{figure}[H]
\centering
\includegraphics[scale=0.95]{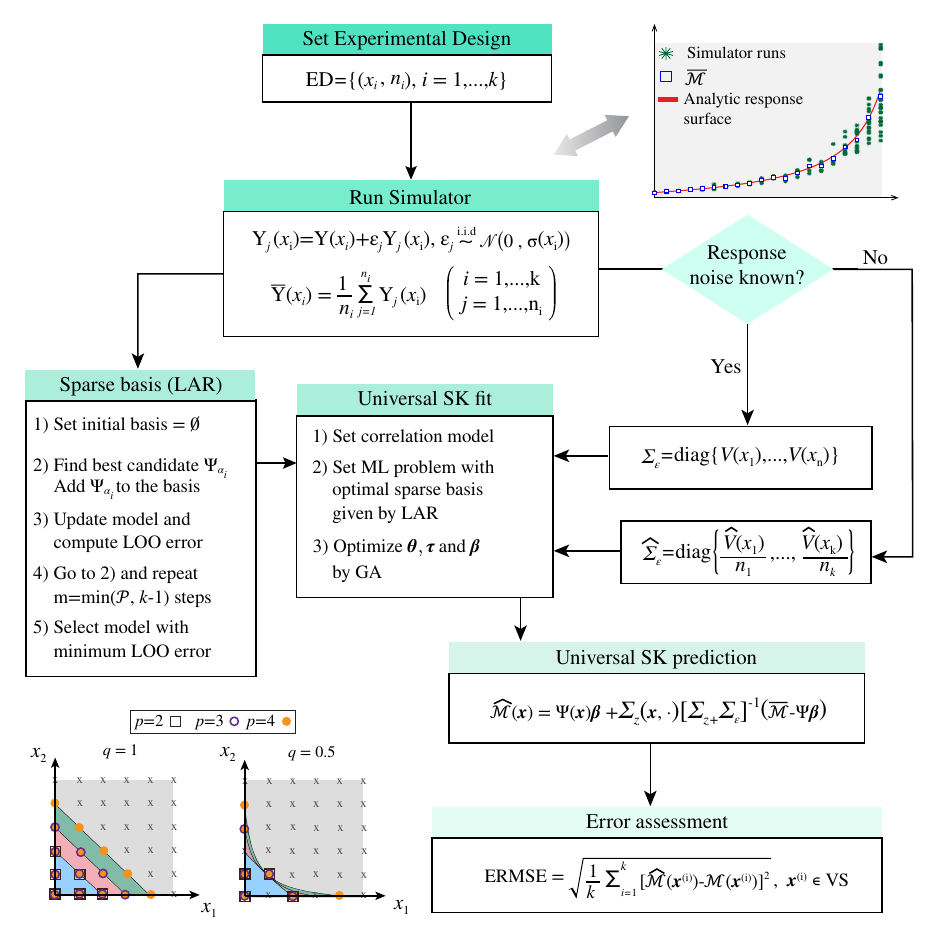}
\caption{Flowchart of the proposed PCE SK surrogate modelling approach for stochastic simulators.}
\label{flowchart}
\end{figure}

 This section presents the proposed PCE SK surrogate modelling approach, synthetically described in Fig. \ref{flowchart}. One of the most remarkable features of the proposed model regards the automatic optimal definition  of the trend  model $\bs{f}$ for Universal SK. This is achieved by the PCE meta-model introduced in Eq.~(\ref{PCE}) using the hyperbolic truncation and the LAR algorithm overviewed in Sections~\ref{PCE_sec}  and \ref{sec_LAR}, respectively. Note that the resulting approach provides  a trend  model with high flexibility in virtue of the LAR algorithm, in such a way that the terms in the PCE expansion are automatically tuned depending on the specific problem at hand. For the development of the proposed meta-model, once the optimal set $\bs{\Psi}$ is identified, a universal SK surrogate model is fitted (see Eq.~(\ref{l-SK})) through a GA global optimization approach  to avoid  stagnation in local extrema.  Then, the prediction of the response surface and the associated MSE are extracted from  Eqs.~(\ref{larpcesk_pred}) and (\ref{larpce_mse_pred}) respectively.

Under the common assumption in the literature that $\left\{ \eps_l \left( \bs{x}^{(i)} \right) \right\}_{l=1}^{n_i}$ are i.i.d. and independent of the random field $\mathcal{Z}$~\cite{Ank-2010,Che-2021,Che-2012}, the proposed optimal LAR-PCE based SK predictor reads:

\begin{equation}\label{pre_PCE SK}
\hat{\mathcal{M}}^{PCE SK} \left( \bs{x} \right)  = {\bs\Psi} \left( \bs{x} \right) {\bs{\beta}}+\bs{\Sigma}_{\mathcal Z} \left( \bs{x},\cdot \right) ^\textrm{T} \left[\bs{\Sigma}_{\mathcal Z}+\bs{\Sigma}_{\eps} \right]^{-1} \left(\overline{{\mathcal{M}}}- \bs{\Psi}{\bs{\beta}} \right),
\end{equation}

\noindent and the corresponding  optimal MSE reads 
\begin{equation}\label{mse-pre_PCE SK}
{\rm MSE} \left( \hat{\mathcal{M}}^{PCE SK} \right) = \bs{\Sigma}_{\mathcal Z} \left( \bs{x},\bs{x} \right) - \bs{\Sigma}_{\mathcal Z} \left( \bs{x},\cdot \right) ^\textrm{T} \left[ \bs{\Sigma}_{\mathcal Z}+\bs{\Sigma}_{\eps} \right]^{-1} \bs{\Sigma}_{\mathcal Z} \left( \bs{x},\cdot \right).
\end{equation} 

It is straightforward to check that the estimator proposed above is the best linear estimator of the form $\hat{\cal M}^{PCE SK} \left( {\bs{x}} \right) = \lambda_0 \left( {\bs{x}} \right)+\bs{\lambda}^\textrm{T} \left( {\bs{x}} \right)\overline{\cal M}$, where $\lambda_0$ and $\bs{\lambda}$ 
are weights depending on $\bs{x}$ and chosen in the sense of minimizing the MSE. Moreover, the estimator is unbiased. This can be deduced by following the general reasoning given in~\cite{Ste-1999} as follows: 

\begin{equation*}
\resizebox{\textwidth}{!}{$
\begin{array}{c}
{\rm{MSE}}={\mathbb{E}}\left[\left({\cal M}^{K}\left( \bs{x} \right)-\hat{\cal M}^{PCE SK}\left( \bs{x} \right)\right)^2 \right]= \mathbb{E}\left[\left({\bs{\Psi}}\left( \bs{x} \right){\bs{\beta}}+{\mathcal Z}\left( \bs{x} \right)-\lambda_0\left( \bs{x} \right)-\bs{\lambda}^\textrm{T}\left( \bs{x} \right) \overline{\mathcal M}\pm\bs{\lambda}^\textrm{T}\left( \bs{x} \right){\bs{\Psi}}{\bs{\beta}}\right)^2\right]
\\\\ 
= \left({\bs{\Psi}}\left( \bs{x} \right){\bs{\beta}}-\lambda_0\left( \bs{x} \right)-\bs{\lambda}^\textrm{T}\left( \bs{x} \right){\bs{\Psi}} {\bs{\beta}}\right)^2 +{\rm{Var}}\left({\mathcal Z}\left( \bs{x} \right)-\bs{\lambda}^\textrm{T}\left( \bs{x} \right)\overline{\cal M}+\bs{\lambda}^\textrm{T}\left( \bs{x} \right){\bs{\Psi}}{\bs{\beta}}\right) + \mathbb{E}\left[\left({\mathcal Z}\left( \bs{x} \right)-\bs{\lambda}^\textrm{T}\left( \bs{x} \right)\overline{\cal M}+\bs{\lambda}^\textrm{T}\left( \bs{x} \right){\bs{\Psi}}{\bs{\beta}}\right)\right]^2
\\\\
=\left({\bs{\Psi}}\left( \bs{x} \right){\bs{\beta}}-\lambda_0\left( \bs{x} \right)-\bs{\lambda}^\textrm{T}\left( \bs{x} \right){\bs{\Psi}}{\bs{\beta}}\right)^2+\Sigma_{\cal Z}\left( \bs{x},\bs{x} \right)+\dis\sum_{i,j=1}^{k}\lambda_i\left( \bs{x} \right)\lambda_j\left( \bs{x} \right){\rm{Cov}} \left({\cal Z}\left( \bs{x}^{(i)} \right),{\cal Z}\left( \bs{x}^{(j)} \right)\right)+\dis\sum_{i,j=1}^{k}\lambda_i\left( \bs{x} \right)\lambda_j\left( \bs{x} \right){\rm{Cov}} \left(\overline{\eps}\left( \bs{x}^{(i)} \right),\overline{\eps}\left( \bs{x}^{(j)} \right)\right)
\\\\
-2\dis\sum_{i=1}^{k}\lambda_i\left( \bs{x} \right)\lambda_j\left( \bs{x} \right){\rm{Cov}} \left({\cal Z}\left( \bs{x}^{(i)} \right),{\cal Z}\left( \bs{x} \right))\right)- 2\dis\sum_{i=1}^{k}\lambda_i\left( \bs{x} \right)\lambda_j\left( \bs{x} \right){\rm{Cov}} \left({\cal Z}\left( \bs{x}^{(i)} \right),\overline{\eps}\left( \bs{x}^{(i)} \right)\right)=
\\\\
=\left({\bs{\Psi}}\left( \bs{x} \right){\bs{\beta}}-\lambda_0\left( \bs{x} \right)-\bs{\lambda}^\textrm{T}\left( \bs{x} \right){\bs{\Psi}}{\bs{\beta}}\right)^2+\Sigma_{\cal Z}\left( \bs{x},\bs{x} \right)+\bs{\lambda}^\textrm{T}\left( \bs{x} \right)\left[\Sigma_{\cal Z}+\Sigma_{\eps}\right]\bs{\lambda}\left( \bs{x} \right)-2\bs{\lambda}^\textrm{T}\left( \bs{x} \right)\Sigma_{\cal Z}({\bs{x}},\cdot),
\end{array}
$}
\end{equation*}

\noindent where $\lambda_i \left( {\bs{x}} \right) $ denotes the $i$-th component of $\bs{\lambda} \left( {\bs{x}} \right) $. Minimising the MSE with respect to $\lambda_0$  and $\bs{\lambda}$, i.e. taking $\partial  \rm{MSE}/\partial \lambda_0 = \partial \rm{MSE}/\partial \bs{\lambda} = 0$,  extract the solutions $\lambda_0 \left( {\bs{x}} \right)  = \Psi  \left( {\bs{x}} \right) {\bs{\beta}}-\bs{\lambda}^\textrm{T} \left( {\bs{x}} \right) {\bs{\Psi}}{\bs{\beta}}$ and $\bs{\lambda} \left( {\bs{x}} \right)  = \left[\Sigma_{\cal Z}+\Sigma_{\eps}\right]^{-1}\Sigma_{\cal Z} \left( {\bs{x}},\cdot \right) $ can be readily extracted. With these weights, the optimal linear estimator and MSE are directly given by Eqs.~(\ref{pre_PCE SK}) and (\ref{mse-pre_PCE SK}), respectively. Furthermore, following that ${\cal Z}$ is a zero-mean stochastic process, i.e.~$\mathbb{E}\left[{\cal Z}\left( \bs{x} \right)\right]=0$, and $\mathbb{E}\left[ \overline{\mathcal{M}}- {\bs{\Psi}}{\bs{\beta}} \right]=0$, it follows that:
\begin{equation}
{\mathbb{E}}\left[\hat{\cal M}^{PCE SK} \left( {\bs{x}} \right) - {\cal M}^K\left( \bs{x} \right)\right] = {\bs\Psi} \left( {\bs{x}} \right) {\bs{\beta}}+\Sigma_{\cal Z} \left( {\bs{x}},\cdot \right)^\textrm{T}\left[\Sigma_{\cal Z}+\Sigma_{\eps}\right]^{-1}\mathbb{E}\left[\overline{\mathcal{M}}- {\bs{\Psi}}{\bs{\beta}}\right] -{\bs{\Psi}} \left( {\bs{x}} \right) {\bs{\beta}}-\mathbb{E}\left[{\cal Z} \left( {\bs{x}} \right) \right]=0,
\end{equation}
\noindent and, therefore, it is extracted that $\hat{\cal M}^{PCE SK}\left( \bs{x} \right)$ is an unbiased estimation for ${\cal M}^{K}\left( \bs{x} \right)$. 
The estimation of the hyper-parameters of the proposed LAR-PCE SK surrogate model requires maximizing the likelihood function described in Eq.~(\ref{l-SK}), where the information matrix $\mathbf{F}$ is simply replaced by the optimal trend matrix $\bs{\Psi}$. However, some considerations must be made. Firstly, note that the covariance matrix of the noise in the response surface (intrinsic bias between the surrogate model and the true simulator's response), $\Sigma_{\varepsilon}$, in Eq.~(\ref{l-SK}) is unknown in general. Nonetheless, $\varepsilon$ is commonly assumed in practice to be normally distributed, so one can write: \begin{equation}\label{esti_vi}    
\left\{\eps_l\left( \bs{x}^{(i)} \right)\right\}_{l=1}^{n_i}\sim {\cal N} \left( 0,V_i \right),\quad \mbox{with}\, \,V_i\,\, \mbox{ estimated as } \quad V_i = \frac{1}{{n_i-1}}\dis\sum_{l=1}^{n_i}\left({\cal M}_l \left( {\bs{x}}^{(i)} \right) - \overline{\cal M}_l \left( {\bs{x}}^{(i)} \right) \right)^2,
\end{equation}
\noindent leading to a positive definite diagonal covariance matrix $\Sigma_{\eps}$ (see Theorem 1 in \cite{Ank-2010}).
Secondly, it should be noted that, when calibrating a Kriging meta-model for deterministic simulators, Eq.~(\ref{lk}) can be simplified as a function of $\bs{\theta}$ by leveraging on the BLUE, while in SK this fact is no longer valid. Instead, ${\bs{\beta}}$ must be first approximated as $\hat{\bm a}$ in Eq.~(\ref{sol}). This  entails limited computational overhead because the trend term has been optimally adjusted beforehand by LAR-PCE. Then, the extrinsic variance ${\sigma}^2$ needs to be simultaneously optimized with the hyper-parameters ${\bs\theta}$ as:
\begin{equation}\label{opt}
\hat{\bs{\theta}},\,\hat{{\sigma}}^2 = \underset{\bs{\theta}, \; {\sigma}^2}{\textrm{arg\,max}}\,{\rm log} \left( {\cal L} \, \left( {\bs{\beta}},{\sigma}^2,\bs{\theta} \right) \right).
\end{equation}
Finally, and based on Eq.~(\ref{opt}), the optimal estimate for ${\bs{\beta}}$ as a function of $\hat{\bs{\theta}}$ and $\hat{\sigma}^2$ can be written similarly to Eq.~(\ref{calibrar}) as:
\begin{equation}\label{opt_betahat}
\hat{\bs{\beta}} \left( \hat{\bs{\theta}},\hat{{\sigma}}^2,{\eps} \right) = \left[ \bs{\Psi}^\textrm{T}\left(\hat{\bs{\Sigma}}_{\cal Z}+\hat{\bs{\Sigma}}_\eps \right)^{-1}\bs{\Psi}\right]^{-1}\bs{\Psi}^\textrm{T}\left(\hat{\bs{\Sigma}}_{\cal Z}+{\hat{\bs{\Sigma}}}_\eps\right)^{-1}\overline{\cal M}.
\end{equation}
 Once the optimal parameters have been derived, the optimal predictor can be re-written as:
\begin{equation}\label{larpcesk_pred}
\hat{\mathcal{M}}^{PCE SK} \left( \bs{x} \right) = {\bs\Psi} \left( \bs{x} \right) \hat{\bs{\beta}} + \hat{\bs{\Sigma}}_{\mathcal Z} \left( \bs{x},\cdot \right)^\textrm{T} \left[ \hat{\bs{\Sigma}}_{\mathcal Z}+{\hat{\bs{\Sigma}}}_{\eps} \right]^{-1} \left( \overline{{\mathcal{M}}}- \bs{\Psi}\hat{\bs{\beta}} \right).
\end{equation}
\noindent with the MSE estimator
\begin{equation}\label{larpce_mse_pred}
\widehat{\rm MSE} \left( \bs{x} \right) = \hat{\bs{\Sigma}}_{\mathcal Z} \left( \bs{x},\bs{x} \right) - \hat{\bs{\Sigma}}_{\mathcal Z} \left( \bs{x},\cdot \right) ^\textrm{T} \left[  \hat{\bs{\Sigma}}_{\mathcal Z}+ \hat{\bs{\Sigma}}_{\eps} \right]^{-1}  \hat{\bs{\Sigma}}_{\mathcal Z} \left( \bs{x},\cdot \right)+{\bs\gamma}^\textrm{T} \left( \bs{\Psi}^\textrm{T}  \left[ \hat{\bs{\Sigma}}_{\mathcal Z}+ \hat{\bs{\Sigma}}_{\eps}\right]^{-1}\bs{\Psi} \right)^{-1}{\bs\gamma}, 
\end{equation}
\noindent   and $\bs{\gamma}= {\bs f} \left( \bs{x} \right) - \bs{\Psi}^\textrm{T} \left[  \hat{\bs{\Sigma}}_{\mathcal Z}+ \hat{\bs{\Sigma}}_{\eps} \right]^{-1}  {\hat{{\bs{\Sigma}}}_{\cal Z} \left( \bs{x},\cdot \right)}$. The expression (\ref{larpce_mse_pred})
follows from the development in \cite{Ste-1999} for standard Kriging by replacing the covariance matrix of the random field by  $\hat{\bs{\Sigma}}_{\mathcal Z}+ \hat{\bs{\Sigma}}_{\eps}$.

Let us lastly focus on the common assumption of the stochastic process ${\cal Z}$ following a Gaussian distribution. In this case, when considering an ED of $k$ samples, the hypotheses in Eq.~(\ref{esti_Z}) imply that $\left({\cal M}\left( \bs{x} \right),\,\overline{\cal M} \left({\bs{x}}^{(1)} \right),\ldots,\,\overline{\cal M} \left( {\bs{x}}^{(k)} \right) \right) \sim {\cal N}_{k+1} \left( \bs{\mu},\bs{\Sigma} \right)$ follow a multivariate Gaussian distribution with mean vector and covariance matrix given by:
\begin{equation}
\bs{\mu} = \left( {\cal M}\left( \bs{x} \right),{\bs\Psi}(\bs{x}) \, {\bs{\beta}} \right) \, \in \, \RR^{k+1}\quad\mbox{and}\quad 
{\bs\Sigma}=\left(
\begin{array}{cc}
{\sigma}^2 &{\bs\Sigma}_{\cal Z}(\bs{x},\cdot)^\textrm{T}\\
{\bs\Sigma}_{\cal Z}(\bs{x},\cdot) & {\bs\Sigma}_{\cal Z}+{{\bs\Sigma}}_{\eps}
\end{array}
\right) \, \in \, \RR^{(k+1) \, \times \, (k+1)},
\end{equation}

\noindent and Eq.~(\ref{pre_PCE SK}) coincides with the conditional mean of ${\cal M}\left( \bs{x} \right)$ given $\overline{\cal M}$, i.e. $\mathbb{E}\left[\left. {\cal M}\left( \bs{x} \right) \right| \overline{\cal M}\right]$ (see~\cite{San-2003}).

To conclude, the quality of the proposed metamodel must be assessed. In practical contexts where the computational burden is high, cross-validation techniques are popular approaches for evaluation the accuracy of stochastic simulators (refer to refernce \cite{Kle-2018} for an extensive discussion). Nevertheless, in this work, to conduct a robust assessment of the proposed PCE SK metamodel, its accuracy has been evaluated on the basis of a dense validation set ${\rm VS}:=\lbrace \bs{x}^{(1)},\ldots,\bs{x}^{(K)}\rbrace$,  with $K \gg k$, covering the parameter space, and using diverse quality metrics reported hereafter.
\section{Numerical results and discussion}\label{numerical}\label{Sect_4}

This section presents  an extensive numerical assessment of the proposed LAR-PCE  SK meta-model. The evaluation includes a detailed comparison against ordinary SK through  various well-known  benchmark examples from the  literature. These presented numerical studies cover both single (Section~\ref{MM1}) and multi-dimensional functions (Sections~\ref{P2_problem_section} and \ref{P3_problem_section}) contaminated by noise. For each case study, given a fixed computational budget, several sampling strategies to obtain the ED are discussed. These approaches range from (i) using a low number of sampling points with many replications to represent the uncertainty to (ii) strategies with numerous sampling points with less budget devoted to the characterization of the uncertainty. Additionally, intermediate strategies that fall between these two extremes are also explored. The purpose of defining different sampling scenarios is to demonstrate that the proposed LAR-PCE SK method consistently outperforms ordinary SK, regardless of whether a sparse and deep ED (with few points but high accuracy) or a denser but shallower ED is chosen. The computational effort allocated at each point is proportional to its corresponding noise variance in every instance. In all the analyses, the model performance is assessed by means of two error metrics, namely the Empirical Root Mean Squared Error (ERMSE) and the Normalized Maximum Absolute Error (NMAE) over a large VS given by: 

\begin{equation}\label{errormetrics}
{\rm ERMSE} = \sqrt{\frac{1}{K} \dis\sum_{i=1}^{K} \left(\widehat{\mathcal{M}} \left( \bs{x}^{(i)} \right) - \mathcal{M} \left( \bs{x}^{(i)} \right) \right)^2}, \quad {\rm NMAE} = \underset{1 \leq i \leq K}{\max} \frac{\left| \widehat{\mathcal{M}} \left( \bs{x}^{(i)} \right) - \mathcal{M} \left( \bs{x}^{(i)} \right) \right|}{K\sigma_{VS}}, \ x^{(i)} \in VS,
\end{equation}

\indent where:

\begin{equation}
\sigma_{VS} = \sqrt{\frac{1}{K} \dis\sum_{i=1}^{K} \left(\overline{\mathcal{M}} \left( \bs{x}^{(i)} \right) - \mathcal{M} \left( \bs{x}^{(i)} \right) \right)^2}.
\end{equation}

Note that the ERMSE metric evaluates the global quality of the meta-model, while the NMAE metric is an indicator of the local error.

\subsection{Case study I: M/M/1 Queue}\label{MM1}

The M/M/1 queue model is a simple yet powerful model used to analyse the behaviour of a wide range of real-world queuing systems~\cite{zheng2000some}, from call centers to manufacturing lines. The model represents a system where there is a single server and a line of entities to be served (e.g., packets, customers, etc.). In this system, entities arrive randomly at a certain rate, called arrival rate, and wait in a queue until the server becomes available. The waiting time is determined by the arrival rate and by the service rate, which refer to the time it takes for the server to process an entity and complete its service, respectively. 

The present case study consists of the simulation of an M/M/1 queue where the service rate is fixed at one and the arrival rate ranges from $0.3$ to $0.9$. The response surface $y$ represents the expected waiting time as a function of the arrival rate $x$. The accuracy of the simulation at each design point is determined by its time run-length $T$, that is, the period of time over which the system's behaviour is observed and analysed. This system is one of the most studied examples in the SK literature~\cite{Ank-2010,Chen-2014,Xie-2020}. It is particularly notable for its distinct and progressively increasing intrinsic variance across the parameter space. This characteristic naturally provides a formidable yet simple example to understand the potential benefits of the proposed approach. Despite the stochastic nature of the model, analytical expressions for the response surface and its intrinsic variance are known~\cite{Ank-2010, Chen-2014}:

\begin{equation}
\label{MM1_eq}    
y(x) = \frac{1}{1-x}, \ V(x) = \frac{2x \left(1+x \right)}{T \left( 1-x \right)^4}.
\end{equation}

In the following, we introduce two variations of the experiment. In the first one in Section~\ref{MM1_One_rep}, the intrinsic variance is assumed as known following Eq.~(\ref{MM1_eq}), and   one single simulation is run for each design point. In the second variation in Section~\ref{MM1_Budget}, the intrinsic variance of the simulator is estimated by performing multiple replications at each of the ED points. For both cases, the results obtained by the classical Ordinary SK are compared with those obtained by the proposed LAR-PCE SK approach.

\subsubsection{Known intrinsic noise} \label{MM1_One_rep}

For this first numerical investigation, we will consider the intrinsic variance as $\Sigma_{\eps} = {\rm{diag}}\left( V\left( \bs{x}^{(i)} \right)\right)_{i=1}^k$. Additionally, $\overline{\mathcal{M}}$ is directly taken as $ \left(\mathcal{M} \left( x^{(i)} \right) \right)_{i=1}^k$ since one single replication at each design point is considered in this first numerical study. It is important to note that, although this situation may not reflect the complexity of most real-world applications, it serves as a useful analysis to shed light on the characteristics of the proposed LAR-PCE SK approach. Three different scenarios with varying number of design points and run-lengths but with constant computational resources have been considered. The first one involves 10 equispaced design points with a run-length of 6000; the second one consists of 30 design points with a run-length of 2000; and the third scenario encompasses 50 design points with a run-length of 1200. In the case of LAR-PCE SK, following the recommendations provided in reference~\cite{Sudret-2018}, the maximum degree of the polynomials is considered as $5$, $10$ and $16$ for the scenarios one to three, respectively. It is worth noting that, given the one-dimensional nature of the problem, the LAR algorithm simply selects the optimal degree for the uni-variate regression model. 

\begin{figure}
\centering
\includegraphics[scale = 1]{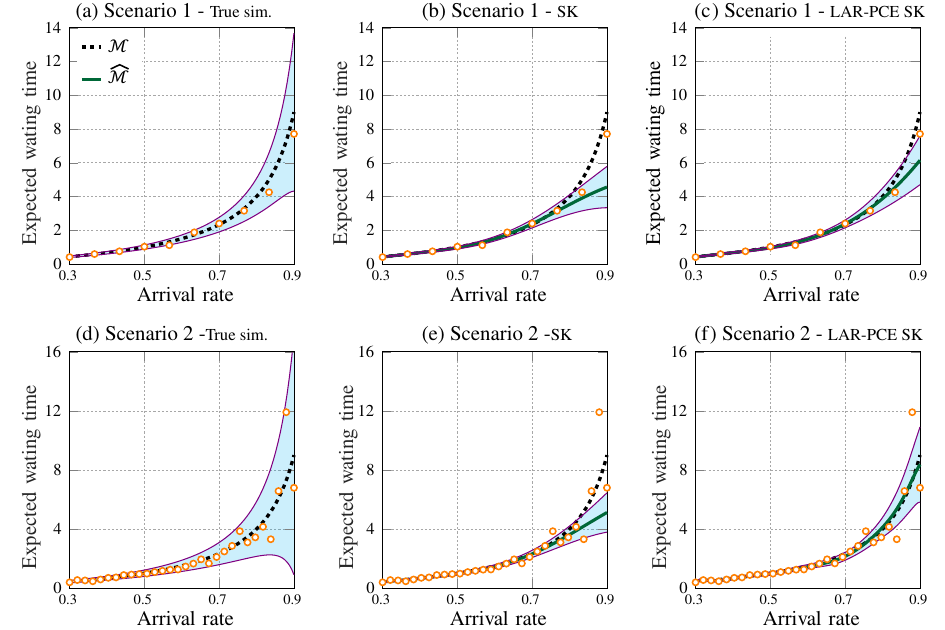}
\caption{Case Study I: M/M/1 example. Comparison between the true simulator (a,d) against the predictions by the proposed LAR-PCE SK (c,f)   and Ordinary SK  (b,e)  under the assumption of known variance (one replication per ED sample) in Scenarios 1 (a,b,c) (10 ED samples and run-length=6000) and 2 (d,e,f) (30 ED samples and run-length=2000). Scatter points denote the ED samples, and light blue shaded areas represent the 95\% confidence interval.}
\label{MM1_results}
\end{figure}

Figures~\ref{MM1_results} (a,d) depict the true response surface along with the analytically computed $95\%$ confidence intervals represented as $\left(y(x)\pm z_{\alpha}\sqrt{V(x)}\right)$ for Scenarios 1 and 2, respectively. Similarly, Figs.~\ref{MM1_results} (b,e) and (c,f) illustrate the estimates produced by Ordinary SK and LAR-PCE SK, respectively. Estimated $95\%$ confidence intervals in Figs.~\ref{MM1_results} (b,e) and (c,f) have been computed as $\left(\hat{\mathcal{M}}^{SK}\pm z_{\alpha}\smash{\sqrt{\widehat{MSE}(x)}}\right)$ and $\left(\hat{\mathcal{M}}^{PCE SK}\pm z_{\alpha}\smash{\sqrt{\widehat{MSE}(x)}}\right)$, respectively. It can be seen that for low values of $x$, where the intrinsic uncertainty is small, both meta-models exhibit a similar behaviour. In contrast, for high values of $x$ where both the QoI and the intrinsic variance experience a rapid growth, Ordinary SK fails to capture the variation in the response, whereas LAR-PCE SK exhibits a significantly better performance. This is to be expected since, as reported in~\cite{Staum-2009}, Ordinary SK virtually ignores those simulation outputs whose estimated intrinsic variance are large compared to the estimated extrinsic variance. Thus, these results demonstrate that introducing a proper trend model effectively mitigates this issue.

\begin{figure}[H]
\centering
\includegraphics[scale = 1]{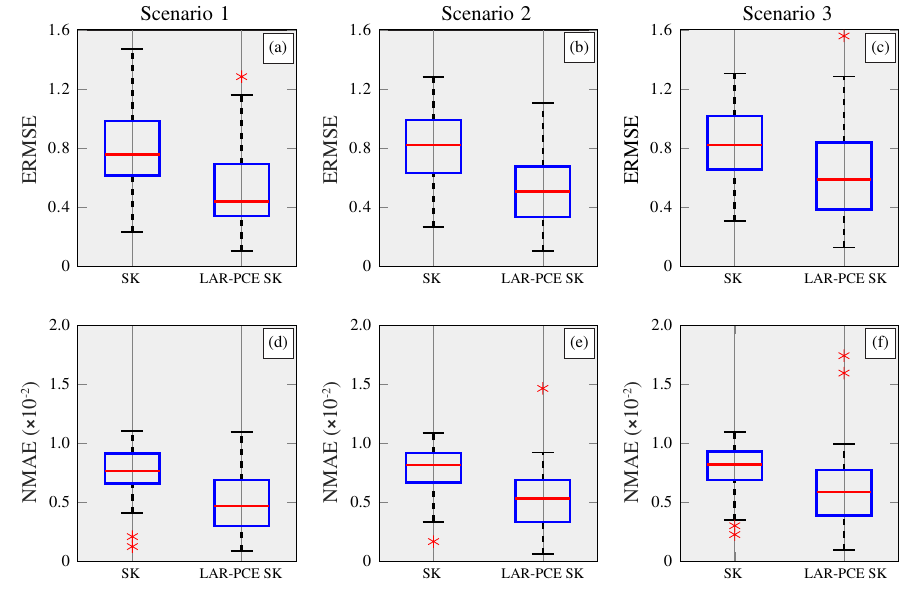}
\caption{Performance analysis of Ordinary SK and LAR-PCE SK in Case Study I in terms of ERMS and NMAE error metrics using 100 known intrinsic noise M/M/1 experiments and sampling Scenarios 1 to 3.}
\label{MM1_results_metrics_one_rep}
\end{figure}

\begin{table}[H]	
\setlength{\tabcolsep}{3pt} 
\newcommand\Tstrut{\rule{0pt}{0,3cm}}         
\newcommand\Bstrut{\rule[-0.15cm]{0pt}{0pt}}   
\centering	
\begin{tabular}{c c c c c c}
\hline
 & \multicolumn{2}{c}{SK} & & \multicolumn{2}{c}{LAR-PCE SK} \\
\cline{2-3} \cline{5-6}
 & Mean ERMSE & Mean NMAE & & Mean ERMSE & Mean NMAE \\
\hline
Scenario 1 & 0.787 & 0.0076 & & 0.516 & 0.0049 \\
Scenario 2 & 0.810 & 0.0079 & & 0.519 & 0.0052 \\
Scenario 3 & 0.832 & 0.0080 & & 0.650 & 0.0074 \\
\hline
\end{tabular}
\footnotesize 
\caption{Summary of the accuracy in the estimates of the M/M/1 model with known intrinsic noise variance by Ordinary SK and LAR-PCE SK adopting sampling Scenarios 1 to 3.}
\label{MM1_metrics_table}
\end{table}
 
To provide a further comparison in performance of LAR-PCE SK with respect to Ordinary SK, the results of the fitting of one hundred independent experiments are reported in Fig.~\ref{MM1_results_metrics_one_rep} and Table~\ref{MM1_metrics_table}. Significant enhancements have been achieved in terms of ERMSE, with improvements ranging from 20\% to 35\% on average. Similarly, the NMAE showed average improvements ranging from 8\% to 35\%. It is noteworthy that, when considering the median values, the improvements are even more significant, ranging from 28\% to 42\% for ERMSE and from 28\% to 38\% for NMAE. Finally, it can be noted in Fig.~\ref{MM1_results} that, in terms of variance estimation, both approaches   similarly   fail to capture the true variability. This issue has been recognized and documented in the SK literature~\cite{Hao-2021, Kle-2016}. Although there are methods available to address this problem (refer e.g. to~\cite{Kle-2016}), they are beyond the scope of our study, in which the primary focus is on the evaluation of the importance of properly selecting the trend model in SK.

\subsubsection{Unknown intrinsic noise}\label{MM1_Budget}

This section introduces a modified version of the previous investigation, where the intrinsic variance is assumed to be unknown. Consequently, following Eq.~(\ref{esti_vi}), the estimation of the $\Sigma_{\varepsilon}$ matrix is carried out through a series of replications at each ED point. In this particular analysis, a total computational budget $C$ of 500 replications to be distributed along the training sites has been considered. Consistently with reference~\cite{Xie-2020}, the number of replications assigned to each design point $x_i$ is determined in a single step as proportional to the model variance following $n_i = \frac{\sqrt{V \left( \bs{x}_i \right)}}{\sum_{i=1}^{k}\sqrt{V \left( \bs{x}_i \right)}}C$. It should be noted that the use of the actual variance to determine the number of replications at each point in our analysis has been chosen for the sake of simplicity. However, alternative methods which employ a multi-step approach for this purpose such as those described in reference~\cite{Chen-2014}   may be considered for bypassing this issue. 

\begin{figure}[H]
\centering
\includegraphics[scale = 1]{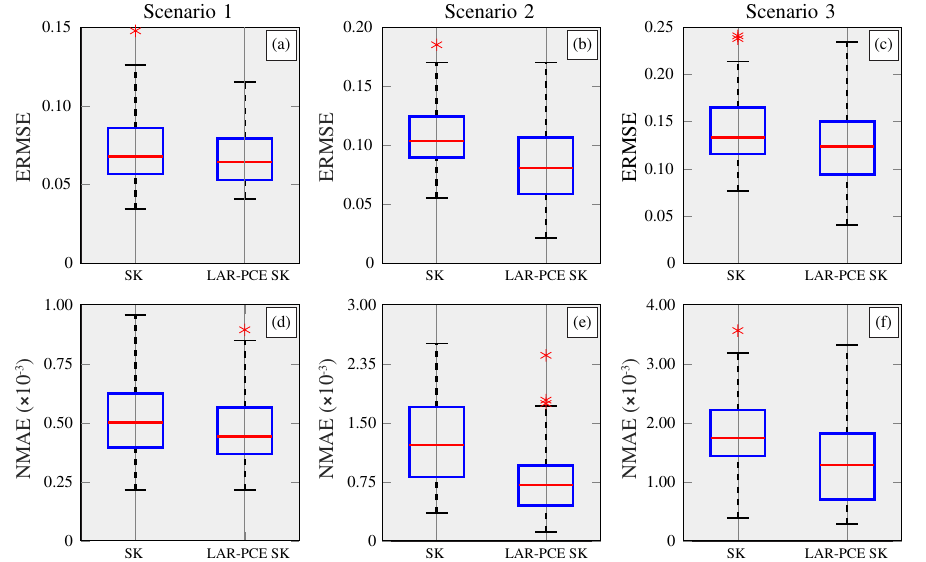}
\caption{Performance analysis of Ordinary SK and LAR-PCE SK in Case Study I in terms of ERMS and NMAE error metrics using 100 unknown intrinsic noise M/M/1 experiments and sampling Scenarios 1 to 3.}
\label{MM1_results_metrics_budget}
\end{figure}

\begin{table}[H]	
\setlength{\tabcolsep}{3pt} 
\newcommand\Tstrut{\rule{0pt}{0,3cm}}         
\newcommand\Bstrut{\rule[-0.15cm]{0pt}{0pt}}   
\centering	
\begin{tabular}{c c c c c c}
\hline
 & \multicolumn{2}{c}{SK} & & \multicolumn{2}{c}{LAR-PCE SK} \Tstrut\Bstrut\\	
\cline{2-3} \cline{5-6}
 & Mean ERMSE & Mean NMAE & & Mean ERMSE & Mean NMAE \Tstrut\Bstrut\\	
\hline
Scenario 1 & 0.073 & 5.16E-04 & & 0.067 & 4.67E-04 \Tstrut\\
Scenario 2 & 0.107 & 1.26E-03 & & 0.084 & 7.53E-04 \\ 
Scenario 3 & 0.142 & 1.80E-03 & & 0.127 & 1.34E-03 \Bstrut\\ 
\hline
\end{tabular}
\footnotesize 
\caption{Summary of the accuracy in the estimates of the M/M/1 model with unknown intrinsic noise variance by Ordinary SK and LAR-PCE SK adopting sampling Scenarios 1 to 3.}
\label{MM1_budget_metrics_table}
\end{table}

The same three scenarios as discussed in Section~\ref{MM1_One_rep} have been considered over a macro-replication of one hundred experiments. The results are summarized in Table~\ref{MM1_budget_metrics_table} and Fig.~\ref{MM1_results_metrics_budget}. As in the previous case, a relevant improvement in performance can be observed in each situation, ranging from 8\% to 21\% on average with respect to ERMSE, and from 9\% to 40\%  with respect to NMAE. Considering the median values, similar conclusions to the previous analyses can be extracted, as  can be seen in Fig.~\ref{MM1_results_metrics_budget}. 

\subsection{Case study II: two dimensional example.}\label{P2_problem_section}

This second case study investigates the egg-box shaped surface suggested in reference~\cite{Xie-2020}. The simulation output at a design point ${\bs{x}} = \left[ x_1,x_2 \right]^\textrm{T}\in \left[-1,1\right]^2$ on the $j$-th replication is given by:

\begin{equation}\label{P2_def}
\mathcal{M}_j \left( x_1,x_2 \right) = \sin \left( 9x_1^2 \right) + \sin \left( 9x_2^2 \right) + \varepsilon_j \left(x_1,x_2 \right), 
\end{equation}

\noindent with $\varepsilon_j \sim {\cal N} \left(0,2+\cos \left( \pi+ \left(x_1+x_2 \right)/2 \right)\right)$.

\begin{figure}[H]
\centering
\includegraphics[scale = 1.3]{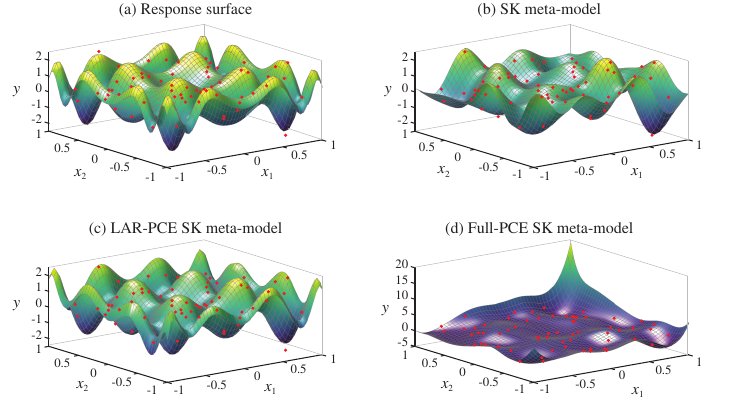}
\caption{Response surface of the Case Study II   (a), and predictions by SK (b), LAR-PCE SK (c) and full PCE SK (d) surrogate models with $k=64$ design points.}
\label{P2_surface}
\end{figure}

Likewise in the previous case study, three different scenarios have been considered   by varying the number of design points but keeping a fixed computational budget $C$ of 1280 replications to be distributed among them. Specifically, the three different settings are characterized by an ED size of $k$=32, 64, and 128 points, respectively. The locations of the design points were determined by adopting the Latin Hypercube Sampling (LHS) algorithm, a commonly used technique for generating diverse and evenly distributed design points in multi-dimensional space. The number of replications $n_i$ at each design point has been adjusted proportionally to its standard deviation as in Section~\ref{MM1_Budget}. For each scenario, the performance of Ordinary SK is benchmarked against the proposed LAR-PCE SK meta-model. Additionally, to stress the importance of choosing an appropriate sparse basis of polynomials in the trend model, a full PCE SK model without LAR has been also included in the comparison. In this case, the entire polynomial basis is included up to a certain degree given by the hyperbolic truncation defined in Section~\ref{PCE_sec}. Regarding the polynomial basis settings, a $q$-norm value of 0.8 has been employed in all cases. With respect to the degree of the polynomial basis to be considered, in compliance with the recommendations by Sudret~\cite{Sudret-2018}, the maximum degree of the polynomial basis for the full PCE SK case has been set at 5, 8, and 10 for Scenarios 1, 2, and 3, correspondingly. On the other hand, due to the sparse nature of the basis selected by the LAR algorithm, it is possible to consider polynomials of slightly higher degree for LAR-PCE SK. Specifically, degrees of 9, 11, and 12 have been considered in Scenarios 1 to 3, respectively. 

\begin{figure}[H]
\centering
\includegraphics[scale = 1.0]{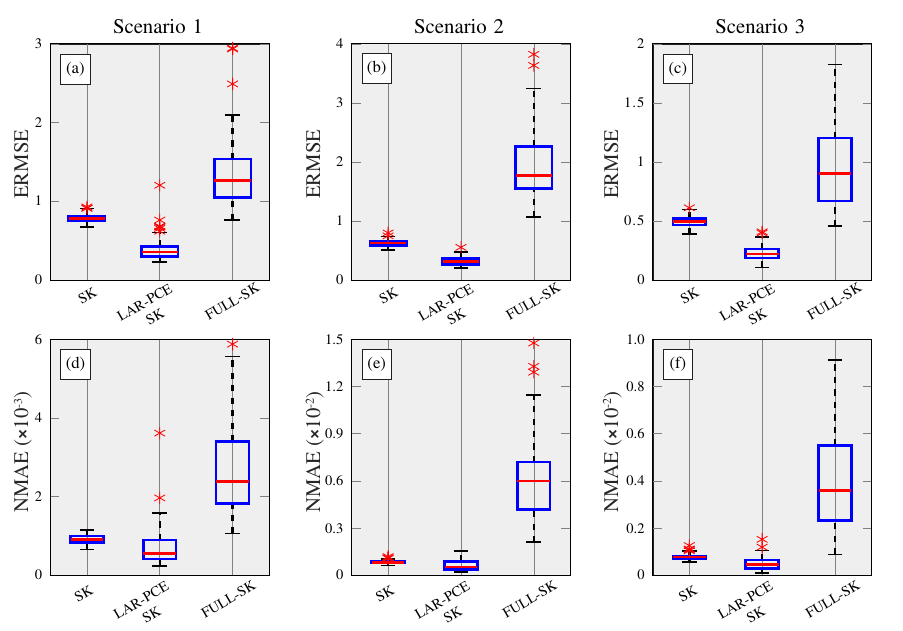}
\caption{Performance analysis of Ordinary SK, LAR-PCE SK and full PCE SK in terms of ERMS and NMAE error metrics using 100 experiments of Case Study II and sampling Scenarios 1 to 3.}
\label{P2_metrics}
\end{figure}

Figure~\ref{P2_surface} illustrates the outcomes of an arbitrary replication of Scenario 2. It depicts the true response surface to be approximated (a), alongside the meta-models constructed using Ordinary SK (b), LAR-PCE SK (c), and full PCE SK (d) with a total of $k=64$ design points. Figure~\ref{P2_metrics} and Table~\ref{P2_metrics_table} present the results of the comparison between the different approaches discussed in the considered scenarios over a macro-replication of one hundred independent experiments. In terms of performance, LAR-PCE SK consistently outperforms the other methods, while full PCE SK showed the poorest performance based on both ERMSE and NMAE metrics. Comparing Ordinary SK with the proposed LAR-PCE SK, average reductions around 50\% in ERMSE are found in all the considered scenarios. With respect to the local accuracy, NMAE enhancements vary from 23\% in Scenarios 1 and 2, to 36\% in Scenario 3. Taking into account the median, similar conclusions can be inferred. On the other hand, it is noteworthy that blindly setting an inadequate polynomial basis as the trend model resulted in notable errors and instability in the results across all the scenarios, highlighting the crucial importance of a proper basis selection. Remarkably, LAR-PCE SK included a modest number of terms in the basis with a median of 9, 11, and 14 terms for Scenarios 1 to 3, respectively. In contrast, the complete expansion included 17, 32, and 49 polynomials in the basis for each respective scenario.

\begin{table}[H]	
\setlength{\tabcolsep}{3pt} 
\newcommand\Tstrut{\rule{0pt}{0,3cm}}         
\newcommand\Bstrut{\rule[-0.15cm]{0pt}{0pt}}   
\centering	
\begin{tabular}{c c c c c c}
\hline
 & \multicolumn{2}{c}{SK} & & \multicolumn{2}{c}{LAR-PCE SK} \Tstrut\Bstrut\\	
\cline{2-3} \cline{5-6}
 & Mean ERMSE & Mean NMAE & & Mean ERMSE & Mean NMAE \Tstrut\Bstrut\\	
\hline
Scenario 1 & 0.782 & 9.13E-04 & & 0.385 & 6.98E-04 \Tstrut\\
Scenario 2 & 0.623 & 8.32E-04 & & 0.316 & 6.42E-04 \\ 
Scenario 3 & 0.495 & 7.69E-04 & & 0.228 & 4.90E-04 \Bstrut\\  
\hline
\end{tabular}
\footnotesize 
\caption{Summary of the accuracy  in   Case Study  II adopting sampling Scenarios 1 to 3.}
\label{P2_metrics_table}
\end{table}

\subsection{Case study III: Ishigami function.}\label{P3_problem_section}

This last case study investigates the Ishigami function~\cite{Schobi-2015}, a classical benchmark function widely used in the literature. This function is characterized by its non-linear and non-monotonic behaviour, making it a challenging test for meta-modelling techniques. The function takes as input a three-dimensional vector $\bs{x} = \left[ x_1,x_2,x_3 \right]^\textrm{T}$ with each component $x_i$ ranging from $-\pi$ to $\pi$, and outputs a scalar value $\mathcal{M}_j\left( \bs{x} \right)$ given by the non-linear sum of three sinusoidal functions. For the purpose of this work, a stochastic version of the Ishigami function is defined by incorporating a normally distributed random noise term $\varepsilon_j\left( \bs{x} \right)$ with zero mean and standard deviation proportional to the absolute value of the function. This results in the following stochastic Ishigami function:

\begin{equation}\label{Ishigami_def}
\begin{array}{c}
\mathcal{M}_j\left( x_1,x_2,x_3 \right) = \sin \left( x_1 \right) + 7\sin^2 \left( x_2 \right) + 0.1 x_3^4 \sin \left( x_1 \right) + \varepsilon_j \left( x_1,x_2,x_3 \right), \ x_1, x_2,x_3 \in \left[ -\pi,\pi \right], 
\\
\varepsilon_j \sim {\cal N} \left(0,\sqrt{\left| \sin \left( x_1 \right) + 7\sin^2 \left( x_2 \right) + 0.1 x_3^4 \sin \left( x_1 \right) \right|}\right).
\end{array}
\end{equation}

The ED has been defined with a total computational budget of $C= 2560$. Three distinct scenarios considering different ED sizes have been investigated: $k=64$ (Scenario 1), $k=128$ (Scenario 2), and $k=256$ (Scenario 3). The allocation of replications has been defined proportionally to the standard deviation, as outlined in Section~\ref{P2_problem_section}. On this basis, Ordinary SK, LAR-PCE SK, and full-PCE SK meta-models have been constructed and compared in terms of accuracy. The calibrations for Scenarios 1 to 3 have been performed with different maximum degrees for the polynomial basis. Specifically, the full-PCE SK meta-models have been calibrated with maximum degrees of 4, 6, and 7, respectively. On the other hand, the LAR-PCE SK meta-models have been calibrated with slightly higher maximum degrees of 6, 8, and 9 for Scenarios 1 to 3, respectively. Additionally, a value of 0.8 is specified for the $q$-norm for all the considered scenarios.

\begin{figure}[H]
\centering
\includegraphics[scale = 0.85]{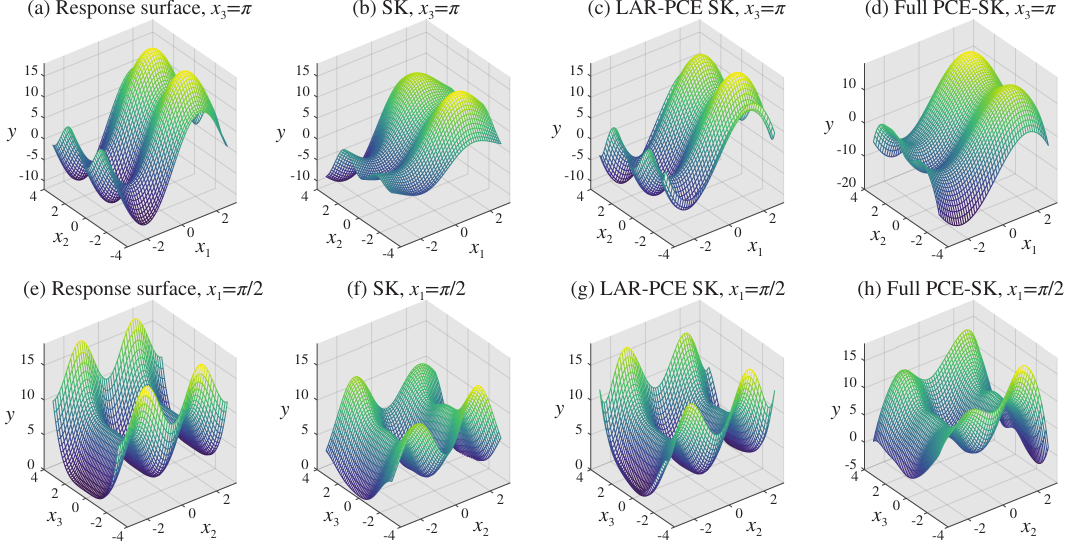}
\caption{Response surfaces of the Ishigami function (a,e), and predictions by Ordinary SK (b,f), LAR-PCE SK (c,g) and full PCE SK (d,h) meta-models with $k = 64$ design points. The top and bottom  row panels  denote the response surfaces obtained for $x_3 = \pi$ and $x_1 = \pi/2$, respectively.}
 \label{Ishigami_results}
\end{figure}

Figure~\ref{Ishigami_results} shows an arbitrary replication of Scenario 1 ($k=64$ design points), displaying the true response surface (a) and the predictions of the meta-models constructed using SK (b), LAR-PCE SK (c), and full PCE SK (d). The comparison among the different approaches based on a replication of 100 independent experiments is presented in Fig.~\ref{P2_metrics} and Table~\ref{P2_metrics_table}. It can be observed that LAR-PCE SK consistently outperforms the other methods in terms of both global (ERMSE) and local (NMAE) accuracy. In this case study, similarly to Case Study II, full PCE SK exhibits the poorest performance. Nonetheless, in this particular case study, the difference between full PCE SK and Ordinary SK is not as pronounced. Instead, when comparing SK with the proposed LAR-PCE SK, significant enhancements are noticeable in terms of ERMSE, with average improvements ranging from 35\% in Scenario 1 to 74\% in Scenario 3. Similar enhancements are also observed when considering the NMAE metric or the improvement in median terms. Lastly, it is worth mentioning the sparsity of the polynomial bases used in the trend model. While the full-PCE expansions consist of a larger number of terms (23, 50, and 68 in Scenarios 1 to 3), the LAR-selected polynomials are significantly more sparse, with only 10, 14, and 17 terms on average, respectively. This allows capturing the trend without encountering over-fitting issues.

\begin{figure}[H]
\centering
\includegraphics[scale = 1.0]{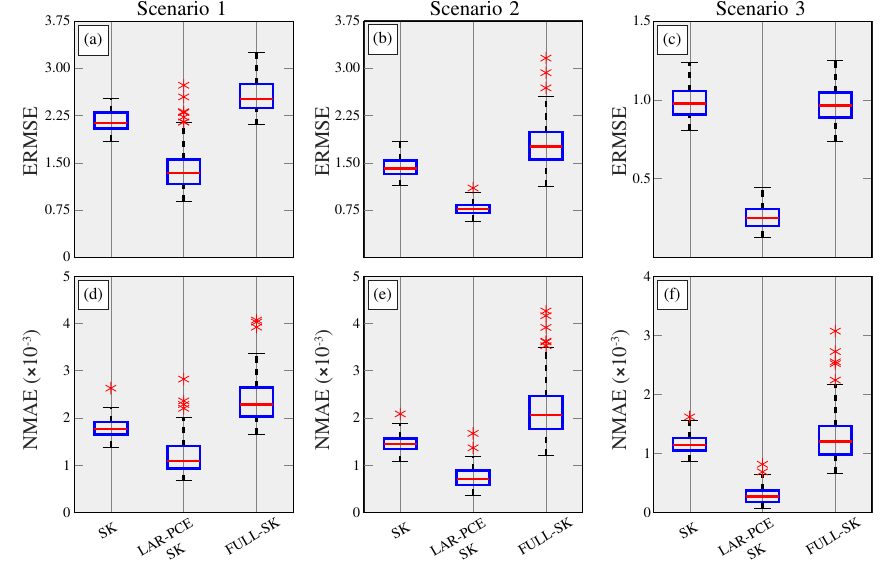}
\caption{Performance analysis of Ordinary SK, LAR-PCE SK and full PCE SK in terms of ERMS and NMAE error metrics using 100 experiments of the Ishigami function and sampling Scenarios 1 to 3.}
\label{Ishigami_metrics_results}
\end{figure}

\begin{table}[H]	
\setlength{\tabcolsep}{3pt} 
\newcommand\Tstrut{\rule{0pt}{0,3cm}}         
\newcommand\Bstrut{\rule[-0.15cm]{0pt}{0pt}}   
\centering	
\begin{tabular}{c c c c c c}
\hline
 & \multicolumn{2}{c}{SK} & & \multicolumn{2}{c}{LAR-PCE SK} \\
\cline{2-3} \cline{5-6}
 & Mean ERMSE & Mean NMAE & & Mean ERMSE & Mean NMAE \\
\hline
Scenario 1 & 2.169 & 1.77E-03 & & 1.418 & 1.20E-03 \\
Scenario 2 & 1.421 & 1.46E-03 & & 0.772 & 7.32E-04 \\ 
Scenario 3 & 0.983 & 1.11E-03 & & 0.257 & 2.92E-04 \Bstrut\\
\hline
\end{tabular}
\footnotesize 
\caption{Summary of the accuracy in the estimates of the Ishigami function for sampling Scenarios 1 ($k=64$), 2 ($k=128$), and 3 ($k=256$).}
\label{Ishigami_metrics_table}
\end{table}

\section{Conclusions}\label{Sect_5}

This study has proposed the use of Universal Stochastic Kriging for constructing surrogate models for stochastic simulators. Specifically, the proposed approach, LAR-PCE SK, adopts adaptive sparse PCE as the trend model in a Universal SK model. With the aim of minimizing over-fitting, the LAR algorithm is adopted to automatically identify the optimal polynomial basis in the PCE, providing the resulting model with local/global prediction abilities and high flexibility for a broad variety of problems. Lastly, the hyper-parameters of the resulting Universal SK model are fitted through a global GA optimization. The effectiveness of the proposed approach has been validated through three classical benchmark case studies, namely: (i) the M/M/1 queue; (ii) a noisy eggbox-shaped surface; and (iii)  a stochastic version of the Ishigami function. The numerical results and discussion have proved the appropriateness of the developed scheme to analyse and develop surrogate models in a wide range of situations, mitigating some of the deficiencies of Ordinary SK when dealing with stochastic processes with high variability. The key findings and contributions of this work include the following:

\begin{itemize}
    \item The proposed LAR-PCE SK approach has demonstrated superior performance  compared to Ordinary SK across various scenarios, including both single and multiple dimensions, with known and unknown intrinsic variance. Benefiting from the  model selection capabilities of the LAR algorithm, LAR-PCE SK has shown particularly promising for multidimensional applications.

    \item Since the computational complexity added by the PCE trend model in Universal SK is marginal, the enhancements in accuracy by the proposed LAR-PCE SK do not entail significant increases in the computational burden compared to Ordinary SK.

    \item The analysis of Case Studies II and III has revealed that the consideration of a full PCE as a trend basis leads to unacceptable results due to over- overfitting. This stresses the need to consider an optimal sparse basis to effectively enhance the performance of the ordinary SK. 
    
\end{itemize}

The presented numerical results and discussion suggest the ability of the proposed methodology to be adapted to a large variety of problems in science and engineering. In this light, future developments regard the extension of the presented formulation to other procedures for sparse basis selection. Another interesting objective for future work involves addressing the inherent underestimation of the variance within the SK methodology. A well-selected trend term could have the potential to improve the ability of the stochastic component in SK to accurately capture the overall stochastic nature of the model. Thus, developing a strategy to effectively handle this underestimation could further enhance the performance and robustness of the surrogate modelling approach.

\section*{Declaration of Competing Interest}
The authors declare that they have no known competing financial interests or personal relationships that could have appeared to influence the work reported in this paper.

\section*{Acknowledgement}

This work has been partially supported by the Spanish Ministry of Science and Innovation [PID2020-116809GB-I00], by the Junta de Extremadura (Spain) through Research Group Grants [GR18023] by the European Regional Development Fund (ERDF) and by ``Junta de Extremadura'' (Ref. IB20040). E. García-Macías was partially supported by the Spanish Ministry of Science and Innovation through the research project ``BRIDGEXT - Life-extension of ageing bridges: Towards a long-term sustainable Structural Health Monitoring'' (Ref. PID2020-116644RB-I00).



\end{document}